\documentclass[11pt, a4paper, twoside]{amsart}

\usepackage[utf8]{inputenc}
\usepackage[english]{babel}
\usepackage[T1]{fontenc}
\usepackage{microtype, fancyhdr, lmodern, xcolor}
\usepackage[a4paper, hmarginratio=1:1]{geometry}

\usepackage{amsfonts,amsmath,amsthm,amscd,amssymb}			
\usepackage{faktor}
\allowdisplaybreaks

\usepackage{url, enumitem}
\usepackage[noadjust]{cite}

\usepackage{graphicx, pdflscape, latexsym}

\input xy
\xyoption{all}

\usepackage{hyperref}

\theoremstyle{plain}
\newtheorem{thm}{Theorem}[section]
\newtheorem{prop}[thm]{Proposition}
\newtheorem{cor}[thm]{Corollary}
\newtheorem{lem}[thm]{Lemma}

\theoremstyle{remark}
\newtheorem{rem}[thm]{Remark}

\theoremstyle{definition}
\newtheorem{defin}[thm]{Definition}

\newcommand\N{\mathbb{N}}
\newcommand\Z{\mathbb{Z}}

\newcommand\R{\mathbb{R}}

\DeclareMathOperator{\Aut}{Aut}	

\newcommand\BB[1]{B_{#1}}	
\newcommand\PB[1]{P_{#1}}	
\newcommand\WB[1]{WB_{#1}}	
\newcommand\VB[1]{VB_{#1}}	
\newcommand\VP[1]{VP_{#1}}	
\newcommand\PC[1]{PC_{#1}}	
\newcommand\UVB[1]{UVB_{#1}}	
\newcommand\UVP[1]{UVP_{#1}}	

\newcommand \F[1]{F_{#1}}

\newcommand\Gr[1]{\langle#1\rangle}		

\newcommand\ie{{\textit{i.e.}}}
\newcommand\ii{i}
\newcommand\inv{^{-1}}

\newcommand\jj{j}

\newcommand\lam[1]{\lambda_{#1}^{\phantom{1}}} 
\newcommand\lamm[2]{\lambda_{#1}^{#2}}	
\newcommand\laminv[1]{\lamm{#1}{-1}}	

\newcommand\nn{n}
\newcommand\nno{{n-1}}

\newcommand\rr[1]{\rho_{#1}^{\phantom{1}}}	
\newcommand\rrq[2]{\rho_{#1}^{#2}}	

\newcommand\sig[1]{\sigma_{#1}^{\phantom{1}}} 
\newcommand\sigg[2]{\sigma_{#1}^{#2}}	
\newcommand\siginv[1]{\sigg{#1}{-1}}	

\begin{document}

\title[Unrestricted virtual braids, fused links and other quotients]{Unrestricted virtual braids, fused links and other quotients of virtual braid groups}

\author[Bardakov]{Valeriy G. Bardakov}
\address{Sobolev Institute of Mathematics, Novosibirsk State University, Novosibirsk 630090, Russia
and Laboratory of Quantum Topology, Chelyabinsk State University, Brat'ev Kashirinykh street 129, Chelyabinsk 454001, Russia;
}
\email{bardakov@math.nsc.ru}

\author[Bellingeri]{Paolo Bellingeri}
\address{Laboratoire de Math\'ematiques Nicolas Oresme, CNRS UMR 6139, Universit\'e de Caen BP 5186,  F-14032 Caen, France.
}
\email{paolo.bellingeri@unicaen.fr}

\author[Damiani]{Celeste Damiani}
\address{Laboratoire de Math\'ematiques Nicolas Oresme, CNRS UMR 6139, Universit\'e de Caen BP 5186,  F-14032 Caen, France.
}
\email{celeste.damiani@unicaen.fr}





\begin{abstract}
We consider the group of unrestricted virtual braids, describe its structure and explore its relations with fused links.
Also, we define the groups of flat virtual braids and virtual Gauss braids and study some of their properties, in particular their linearity.
\end{abstract}

\maketitle

\section{Introduction}

\begin{figure}[htb]
	\centering
		\includegraphics[scale=0.7]{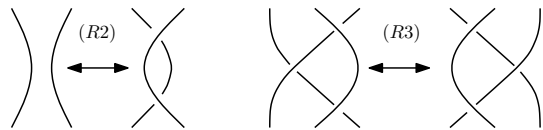} 
	\caption{Classical Reidemeister moves.}
	\label{F:classical}
\end{figure}

\begin{figure}[htb]
	\centering
		\includegraphics[scale=0.7]{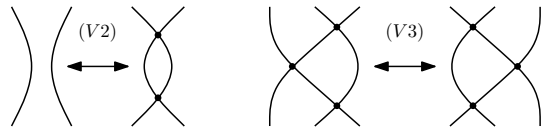} 
	\caption{Virtual Reidemeister moves.}
	\label{F:virtual}
\end{figure}

\begin{figure}[htb]
	\centering
		\includegraphics[scale=0.7]{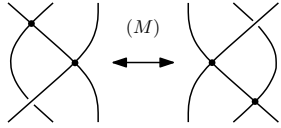} 
	\caption{Mixed Reidemeister move.}

	\label{F:mixed}
\end{figure}

\begin{figure}[htb]
	\centering
		\includegraphics[scale=0.7]{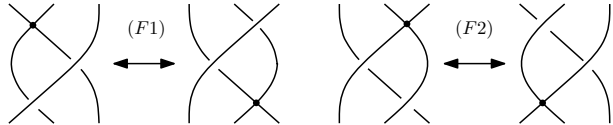} 
	\caption{Forbidden moves of type \eqref{E:F1} (on the left) and type \eqref{E:F2} (on the right).}
	\label{F:forbidden}
\end{figure}

\emph{Fused links} were defined by L.~H.~Kauffman and S.~Lambropoulou in~\cite{KL1}. Afterwards, the same authors introduced their ``braided'' counterpart,
the \emph{unrestricted virtual braids}, and extended S.~Kamada's work~(\cite{Kam}) by presenting a version of Alexander and Markov theorems for these objects~\cite{KaL}.
In the \emph{group of unrestricted virtual braids}, which will be denoted by~$\UVB\nn$, we consider braid-like diagrams in which we allow two kinds of crossing (classical and virtual), and where the equivalence relation is given by ambient isotopy and by the following transformations: classical Reidemeister moves (Figure~\ref{F:classical}), virtual Reidemeister moves (Figure~\ref{F:virtual}), a mixed Reidemeister move (Figure~\ref{F:mixed}), and two moves of type Reidemeister III with two classical crossings and one virtual crossing (Figure~\ref{F:forbidden}). These two last moves are called \emph{forbidden moves}. 

The group~$\UVB\nn$ appears also in~\cite{KMRW}, where it is called \emph{symmetric loop braid group}, being a quotient of the \emph{loop braid group}~$LB_\nn$ studied in~\cite{BWC}, also known as the 
welded braid group~$\WB\nn$.

It has been shown that all fused knots are equivalent to the unknot~(\cite{Kan, Nel}). Moreover, S.~Nelson's proof in~\cite{Nel} of the fact that every virtual knot unknots, when allowing forbidden moves, which is carried on using Gauss diagrams, can be adapted verbatim to links with several components. So, every fused link diagram is fused isotopic to a link diagram where the only crossings (classical or virtual) are the ones involving different components.

On the other hand, there are non-trivial fused links and their classification is not (completely) trivial~(\cite{FK1}): in particular in~\cite{FK2}, A.~Fish and E.~Keyman proved that fused links that have only classical crossings are characterized by their (classical) linking numbers. However, this result does not generalize to links with virtual crossings: in fact it is easy to find non-equivalent fused links with the same (classical) linking number (see Remark~\ref{R:linking}). This answers a question from~\cite[Remark~1]{FK2}, where Fish and Keyman ask whether the classical linking number is a complete invariant for fused links.

The first aim of this note is to give a short survey on above knotted objects, describe unrestricted virtual braids and compare
more or less known invariants for fused links.
In Section~\ref{S:Unrestricted} we give a description of the structure of the group of unrestricted virtual braids~$\UVB\nn$ (Theorems~\ref{T:uvbn} and~\ref{T:theorem1}), answering a question of Kauffman and Lambropoulou from~\cite{KaL}. In Section~\ref{S:Fused} we provide an application of  Theorem~\ref{T:theorem1} showing that any fused link admits
as a representative the closure of a \emph{pure} unrestricted virtual braid (Theorem~\ref{T:fused}); as a corollary we deduce an easy proof of  the theorem of Fish and Keyman cited in previous paragraph. 
In Section~\ref{S:fusedgroup} we construct a representation for $\UVB\nn$ in~$\Aut(N_n)$, the group of automorphisms of the free $2$-step nilpotent group of rank $\nn$ (Proposition~\ref{P:representation}). Using this representation we define a notion of group of fused links and we compare this invariant to other known invariants (Proposition~\ref{P:Group} and Remark~\ref{R:calculation}). 

Finally, in Section~\ref{S:quotients} we describe the structure of other quotients of virtual braid groups: the flat virtual braid group~(Proposition~\ref{P:VirtualFlat} and Theorem~\ref{T:flatH}),
the flat welded braid group (Proposition~\ref{P:fwbn}) and the virtual Gauss braid group (Theorem~\ref{T:gaussH}). As a corollary we prove that 
 flat virtual braid groups and virtual Gauss braid groups are linear and that they have solvable
word problem (the fact that unrestricted virtual braid groups are linear and have solvable
word problem is a trivial consequence of Theorem~\ref{T:theorem1}).

\medskip

\noindent \textbf{Acknowledgments.} The research of the first author was partially supported by Laboratory of Quantum Topology of Chelyabinsk State University (Russian Federation government grant 14.Z50.31.0020),
 RFBR-14-01-00014, RFBR-15-01-00745 and Indo-Russian RFBR-13-01-92697.
The research of the second author was partially supported by French grant ANR-11-JS01-002-01. This paper was started when the first author was in Caen. He thanks the members of the Laboratory of Mathematics of the University of Caen for their invitation and hospitality.
The authors are deeply indebted to the anonymous referee for valuable remarks and suggestions.


\section{Unrestricted virtual braid groups}
\label{S:Unrestricted}

In this Section, in order to define unrestricted virtual braid groups, we will first introduce virtual and welded braid groups by simply recalling their group presentation; 
for other definitions, more intrinsic, see for instance~\cite{Ver01, Bn,Ci,Kam}
 for the virtual case and~\cite{BH,FRR,Kam} for the welded one.

\begin{defin}
The \emph{virtual braid group} $VB_n$ is the group defined by the group presentation
\[ \Gr{ \ \{\sig{\ii}, \rr\ii \mid \ii = 1, \dots, \nno\} \ \vert \ R \ } \]

where $R$ is the set of relations:
\begin{align}
\label{E:R1}\tag{R1} \sig{\ii}  \sig{i+1} \sig{\ii} &= \sig{i+1}\sig{\ii} \sig{i+1}\mbox{,} \quad &\mbox{for } \ii=1,\dots,\nn-2 \mbox{;} \\
\label{E:R2}\tag{R2} \sig{\ii}  \sig{\jj} &= \sig{\jj} \sig{\ii} \mbox{,} \quad &\mbox{for } \vert \ii - \jj \vert \geq 2 \mbox{;} \\
\label{E:R3}\tag{R3} \rr\ii \rr{i+1} \rr\ii &= \rr{i+1} \rr\ii\rr{i+1}\mbox{,} \quad &\mbox{for } \ii=1,\dots,\nn-2 \mbox{;} \\
\label{E:R4}\tag{R4} \rr\ii \rr\jj &= \rr\jj \rr\ii \mbox{,} \quad &\mbox{for } \vert \ii - \jj \vert \geq 2 \mbox{;} \\
\label{E:R5}\tag{R5} \rrq\ii2 &=1\mbox{,} \quad &\mbox{for } \ii=1,\dots,\nno \mbox{;} \\ 
\label{E:R6}\tag{R6} \sig{i}  \rr{j} &= \rr{j} \sig{i}\mbox{,} \quad &\mbox{for } \vert \ii - \jj \vert \geq 2 \mbox{;} \\
\label{E:M}\tag{M} \rr{i}  \rr{i+1} \sig{i} &= \sig{i+1}  \rr{i}   \rr{i+1} \mbox{,} \quad &\mbox{for } \ii=1,\dots,\nn-2.
\end{align}
\end{defin}

We define the \emph{virtual pure braid group}, denoted by~$VP_\nn$, to be the kernel of the map $VB_\nn \to S_\nn$ sending, for every $\ii = 1, 2, \dots, n-1 $, generators $\sig{\ii}$ and $\rr\ii$ to~$(i, i+1)$. A presentation for $VP_\nn$ is given in~\cite{Bar-0}; it will be recalled in the proof of Theorem \ref{T:theorem1} and Proposition~\ref{P:VirtualFlat}.

The welded braid group $WB_n$ can be defined as a quotient of $VB_n$ by the normal subgroup generated by relations
\begin{equation}
\label{E:F1}
\tag{F1} \rr\ii \sig{\ii+1} \sig\ii =  \sig{\ii+1} \sig\ii \rr{\ii+1} \mbox{,} \quad \mbox{for } \ii = 1, \dots, \nn-2.
\end{equation}

\begin{rem}
We will see in Section~\ref{S:Fused} that the symmetrical relations
\begin{equation}
\label{E:F2}
\tag{F2} \rr{\ii+1} \sig{\ii} \sig{\ii+1} =  \sig\ii \sig{\ii+1}  \rr{\ii} \mbox{,} \quad \mbox{for }  \ii = 1, \dots, \nn-2
\end{equation}
do not hold in~$\WB\nn$. This justifies Definition~\ref{D:Unrestricted}.
\end{rem}

\begin{defin}
\label{D:Unrestricted}
We define the \emph{group of unrestricted virtual braids} $\UVB\nn$ as the group defined by the group presentation
\[ \Gr{ \ \{\sig{\ii}, \rr\ii \mid \ii = 1, \dots, \nno\} \ \vert \ R' \ } \]

where $R'$ is the set of relations \eqref{E:R1}, \eqref{E:R2}, \eqref{E:R3}, \eqref{E:R4}, \eqref{E:R5}, \eqref{E:R6}, \eqref{E:M}, \eqref{E:F1}, \eqref{E:F2}.
\end{defin}

The main result of this section is to prove that $\UVB\nn$ can be  described as semi-direct product of a right-angled Artin group and the symmetric group~$S_\nn$:
this way we answer a question posed in \cite{KaL} about the (non-trivial) structure of~$\UVB\nn$.

\begin{thm}
\label{T:uvbn}
Let $X_\nn$ be the right-angled Artin group generated by $x_{i,j}$ for $1\le i\not=j \le \nn$ where all generators commute
except the pairs $x_{i,j}$ and $x_{j,i}$ for~$1 \le i \neq j \le \nn$.
The group $UVB_\nn$ is isomorphic to $X_\nn \rtimes S_\nn$ where $S_\nn$ acts by permutation on the indices of generators
of~$X_\nn$.
\end{thm}

Let
$\nu \colon UVB_\nn \to S_\nn$
be the map defined as follows:
\[
\nu(\sig{\ii}) = \nu(\rr\ii) = (i,  i+1), \quad \mbox{for } i = 1, 2, \dots, \nno.
\]
We will call the kernel of $\nu$ 
\emph{unrestricted virtual pure braid group} and we will denote it by~$\UVP\nn$.
Since $\nu$ admits a natural section, we have that~$\UVB\nn = \UVP\nn \rtimes S_\nn$.

\begin{figure}[htb]
	\centering
		\includegraphics[scale=0.7]{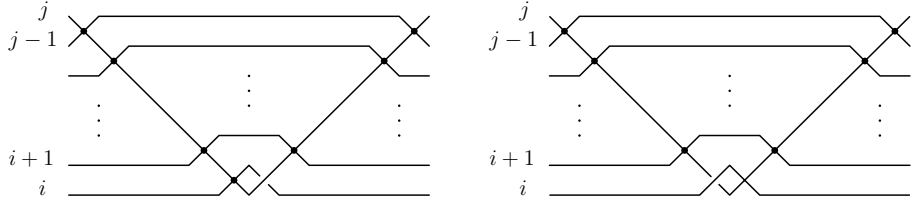} 
	\caption{Elements $\lambda_{\ii, \jj}$ on the right and $\lambda_{\jj, \ii}$ on the left. Here we adopt the convention of drawing braids from left to right.}
	\label{F:lambdaij}
\end{figure}

Let us define some elements of~$\UVP\nn$ (see Figure~\ref{F:lambdaij}). For $\ii = 1, \dots, \nno$:
\begin{equation}
\label{E:lambda}
\begin{split}
&\lam{i,i+1} = \rr\ii \siginv\ii,\\
&\lam{i+1,i} = \rr\ii \lam{i,i+1} \rr\ii = \siginv\ii \rr\ii.
\end{split}
\end{equation}

For $1 \leq \ii < \jj-1 \leq \nno$:
\begin{equation}
\label{E:lambda2}
\begin{split}
&\lam{i,j} = \rr{j-1} \rr{j-2} \ldots \rr{i+1} \lam{i,i+1} \rr{i+1} \ldots \rr{j-2} \rr{j-1}, \\
&\lam{j,i} = \rr{j-1} \rr{j-2} \ldots \rr{i+1} \lam{i+1,i}  \rr{i+1} \ldots \rr{j-2} \rr{j-1}.
\end{split}
\end{equation}

The next lemma was proved in \cite{Bar-0} for the corresponding elements in $\VB\nn$ and therefore is also true in the quotient~$\UVB\nn$.

\begin{lem} 
\label{L:lemma1}
The following conjugating rule is fulfilled in~$\UVB\nn$: for all $1 \leq \ii \neq \jj \leq \nn$ and~$s \in S_n$, 
\[
\iota(s) \lam{\ii ,\jj} \iota(s)\inv = \lam{s(\ii), s(\jj)}
\]
where $\iota \colon S_n \to \UVB\nn$ is the natural section of the map $\nu$ defined in Theorem~\ref{T:uvbn}.
\end{lem}

\begin{cor}
\label{C:cor1}
The group $S_n$ acts by conjugation on the set
$\{ \lam{k,l} ~ \vert 1 ~ \leq k \neq l \leq n \}.$ This action is transitive.
\end{cor}

We prove that the group generated by $\{ \lam{k,l} \ \vert \ 1 \leq k \neq l \leq n \} $ coincides with~$UVP_n$, and then we will find the defining relations. This will show that $\UVP\nn$ is a right-angled Artin group.

\begin{thm} 
\label{T:theorem1}
The group $\UVP\nn$ admits a presentation with generators~$\lam{k,l}$ for $1 \leq k \neq l \leq
n$,
and defining relations:
$\lam{i,j}$ commutes with $\lam{k,l}$ if and only if~$k \not= j$ or~$l \not= i$ .
\end{thm}

 \proof
Since $\UVP\nn$ is a finite index subgroup of $\UVB\nn$  one can apply Reidemeister--Schreier method (see, for example,~\cite[Ch. 2.2]{MKS}) and check that the given set of relations is complete. Remark that most of the relations were already proven in this way in~\cite{Bar-0} for the case of the virtual pure braid group~$\VP\nn$.
 
An easier approach is provided by the following commutative diagram:
\[
\begin{CD}
@. 1 @. 1 \\
@. @VVV  @VVV \\
@. \ker \pi_{\vert \VP\nn} @>>> \ker \pi \\
@. @VVV  @VVV @. \\
1 @>>> \VP\nn @>>>  \VB\nn @>>> S_\nn @>>>1 \\
   @.  @VV{ \pi_{\vert \VP\nn}}V         @VV{\pi}V        @|  \\
1 @>>> \UVP\nn @>>>  \UVB\nn @>>> S_\nn @>>>1 \\
@. @VVV  @VVV \\
@. 1 @. 1 \\
\end{CD}
\] 
where $\pi$ is the canonical projection of $VB_n$ onto $UVB_n$ and $\pi_{\vert \VP\nn}$ its restriction to~$VP_n$.
By definition $ \ker \pi$ is normally generated by elements $\sig\ii \sig\jj \rr\ii  \siginv\jj \siginv\ii \rr\jj$ for $\vert i-j \vert =1$
(we will write $ \ker \pi= \ll \sig\ii \sig\jj \rr\ii  \siginv\jj \siginv\ii \rr\jj \vert  \; \mbox{for} \; \vert i-j \vert =1 \; \gg$).  Since  $ \sig\ii \sig\jj \rr\ii  \siginv\jj \siginv\ii \rr\jj$  belongs to $\VP\nn$ and that
$\VP\nn$ is normal in $\VB\nn$,  we deduce that $ \ker \pi_{\VP\nn}$ coincides with $\ker \pi$.

We recall that, according to~\cite{Bar-0}, $\VP\nn$ is generated by elements $\lam{i,j}$ defined  as follows:
\begin{equation}
\label{E:lambdavirtual1}
\begin{split}
&\lam{i,i+1} = \rr\ii \siginv\ii,\\
&\lam{i+1,i} = \rr\ii \lam{i,i+1} \rr\ii = \siginv\ii \rr\ii.
\end{split}
\end{equation}

For $1 \leq \ii < \jj-1 \leq \nno$:
\begin{equation}
\label{E:lambdavirtual2}
\begin{split}
&\lam{i,j} = \rr{j-1} \rr{j-2} \ldots \rr{i+1} \lam{i,i+1} \rr{i+1} \ldots \rr{j-2} \rr{j-1}, \\
&\lam{j,i} = \rr{j-1} \rr{j-2} \ldots \rr{i+1} \lam{i+1,i}  \rr{i+1} \ldots \rr{j-2} \rr{j-1}.
\end{split}
\end{equation}

and has the following set of defining relations:

\begin{align}
\label{E:RS1}\tag{RS1} \qquad  \lam{i,j}   \lam{k,l} &= \lam{k,l}   \lam{i,j} \\
\label{E:RS2}\tag{RS2} \qquad \lam{k,i}  (\lam{k,j} \lam{i,j}) &= (\lam{i,j} \lam{k,j})   \lam{k,i}.
\end{align}

Moreover, as $\UVB\nn$, $\VB\nn$ can be seen as a semidirect product $\VP\nn\rtimes S_n$, where the symmetric group $S_n$ acts by permutations of indices on $\lam{i,j}$'s~(Lemma~\ref{L:lemma1}).

One can easily verify that relators of type~\eqref{E:F1},  \ie, $\rr{i} \sig{i+1}   \sig{i}   \rr{i+1}   \siginv{i}
\siginv{i+1}$, can be rewritten as:
 \[
 (\rr\ii  \, \laminv{i+1,i+2} \, \rr\ii) 
(\rr\ii \,  \rr{i+1} \,  \laminv{i,i+1} \,  \rr{i+1} \, \rr\ii)
(\rr{i+1} \, \lam{i,i+1} \, \rr{i+1}) 
\lam{i+1,i+2}
\]

and using the conjugating rule given above,
we get, for $i=1,\dots, n-2$,
\[\rr{i} \sig{i+1}  \sig{i}   \rr{i+1}   \siginv{i}  \siginv{i+1}=\laminv{i,i+2} \,  \laminv{i+1,i+2}   \, \lam{i,i+2} \,  \lam{i+1,i+2}.\]

On the other hand one can similarly check that relators of type~\eqref{E:F2}, which are of the form $\rr{i+1}  \sig{i}  \sig{i+1}\rr{i}  \siginv{i+1} 
\siginv{i}$, can be rewritten as $\laminv{i,i+1}   \laminv{i,i+2}   \lam{i,i+1}  \lam{i,i+2}$.

From this facts and from above description of $\VB\nn$ as semidirect product $\VP\nn \rtimes S_n$,
 it follows that any generator of $\ker_{\VP\nn}$
is of the form $g [ \lam{i,j}, \lam{k,j} ] g^{-1} $ or $g[ \lam{i,j}, \lam{i,k} ] g^{-1} $
for $g \in \VP\nn$ and $i,j,k$ distinct.
The group $\UVP\nn$ has therefore the following complete set of relations
\begin{align}
\tag{\ref{E:RS1}} \qquad  \lam{i,j}   \lam{k,l} &= \lam{k,l}   \lam{i,j} \\
\tag{\ref{E:RS2}}  \qquad \lam{k,i}  (\lam{k,j} \lam{i,j}) &= (\lam{i,j} \lam{k,j})   \lam{k,i} \\
\label{E:RS3}\tag{RS3} \qquad \lam{i,j}  \lam{k,j}& =   \lam{k,j} \lam{i,j} \\
\label{E:RS4}\tag{RS4} \qquad  \lam{i,j} \lam{i,k}&= \lam{i,k} \lam{i,j} .\\
\intertext{Using \eqref{E:RS3} and  \eqref{E:RS4}  
we can rewrite relation \eqref{E:RS2} in the form}
\lam{k,j} (\lam{k,i} \lam{i,j}) &= \lam{k,j} (\lam{i,j} \lam{k,i}).\\
\intertext{After cancelation we have that we can replace relation \eqref{E:RS2}
with}
\label{E:RS5}\tag{RS5} \qquad  \lam{k,i} \lam{i,j} &= \lam{i,j}\lam{k,i}
\end{align}
 
This completes the proof.

\endproof

\proof[Proof of Theorem~\ref{T:uvbn}]
The group $X_n$ is evidently isomorphic to $\UVP\nn$ (sending any $x_{i,j}$ into the corresponding~$\lam{i,j}$).
Recall that $\UVP\nn$ is the kernel of the map $\nu \colon \UVB\nn \to S_n$ defined as~$\nu(\sig{\ii}) = \nu(\rr\ii) = (i, i+1)$ for~$\ii = 1, \dots, \nno$.
Recall also that $\nu$ has a natural section~$\iota \colon S_n \to \UVB\nn$, defined as $\iota \big((i, i+1)\big) = \rr\ii$ for~$\ii = 1, \dots, \nno$.
Therefore $\UVB\nn$ is isomorphic to
$\UVP\nn \rtimes S_n$ where $S_n$ acts by permutation on the indices of generators
of $\UVP\nn$ (see Corollary \ref{C:cor1}).

\endproof

We recall that the pure braid group $\PB\nn$ is the kernel of the homomorphism from $\BB\nn$ to the symmetric group $S_n$ sending every generator $\sig\ii$ to the permutation $(i, i+1)$. It is generated by the set $\{a_{ij} \mid 1 \leq i < j \leq \nn\}$, where
\begin{align*}
a_{\ii,\ii+1} &= \sigg\ii{2}, \\
a_{\ii, \jj} &= \sig{\jj-1} \sig{\jj-2} \cdots \sig{\ii+1} \sigg\ii{2} \siginv{\ii+1} \cdots \siginv{\jj-2} \siginv{\jj-1} \mbox{,} \quad  \mbox{for } \ii+1 < \jj \leq \nn.
\end{align*}

\begin{cor}
\label{C:purecommute}
Let $p \colon \PB\nn \to \UVP\nn$ be the canonical map of the pure braid group $\PB\nn$ in~$\UVP\nn$. Then $p(\PB\nn)$ is isomorphic to the abelianization of~$\PB\nn$. 
\end{cor}

 \proof
As remarked in~(\cite[page~6]{Bar-0}), generators $a_{\ii, \jj}$ of $\PB\nn$ can be rewritten in $\VP\nn$ as
\begin{align*}
a_{\ii, \ii+1} &= \laminv{\ii, \ii+1} \laminv{\ii+1, \ii} \mbox{,} \quad \mbox{for } \ii=1, \dots, \nno \mbox{,} \\
a_{\ii, \jj} &= \laminv{\jj-1, \jj} \laminv{\jj-2, \jj} \cdots \laminv{\ii+1, \jj} ( \laminv{\ii, \jj}  \laminv{\jj, \ii}) \lam{\ii+1, \jj} \cdots \lam{\jj-2, \jj}  \lam{\jj-1, \jj} \mbox{,} \quad  \mbox{for } 2 \leq \ii+1 < \jj \leq \nn,
\intertext{and therefore in $\UVP\nn$ we have:}
p(a_{\ii, \ii+1}) &= \laminv{\ii, \ii+1} \laminv{\ii+1, \ii} \mbox{,} \quad \mbox{for } \ii=1, \dots, \nno \mbox{,} \\
p(a_{\ii, \jj})&= \laminv{\jj-1, \jj} \laminv{\jj-2, \jj} \cdots \laminv{\ii+1, \jj} ( \laminv{\ii, \jj}  \laminv{\jj, \ii}) \lam{\ii+1, \jj} \cdots \lam{\jj-2, \jj}  \lam{\jj-1, \jj} \mbox{,} \quad  \mbox{for } 2 \leq \ii+1 < \jj \leq \nn.
\end{align*}

According to Theorem~\ref{T:theorem1}, $\UVP\nn$ is the cartesian product of the free groups of rank $2$  $\F{i, j}= \langle \lam{i,j}, \lam{j,i} \rangle $ for $1 \leq i < j \leq n$. 

For every generator $a_{i, j}$ for $1 \leq i < j \leq n$ of $\PB\nn$ we have that its image is in $\F{i, j}$ and it is not trivial. In fact, $p(a_{i, j}) =  \laminv{\ii, \jj} \laminv{\jj, \ii}$. So $p(\PB\nn)$ is isomorphic to $\Z^{n(n-1)/2}$. The statement therefore follows readily since the abelianized of $\PB\nn$ is $\Z^{n(n-1)/2}$.

 \endproof


\section{Unrestricted virtual braids and fused links}\label{S:Fused}

\begin{figure}[htb]
	\centering
		\includegraphics[scale=0.7]{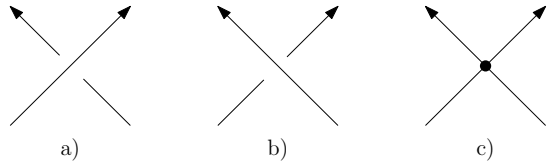} 
	\caption{a) Positive crossing, b) Negative crossing, c) Virtual crossing.}
	\label{F:Crossings}
\end{figure}

\begin{defin}
A \emph{virtual link diagram} is a closed oriented 1-manifold $D$ immersed in~$\R^2$ such that all multiple points are transverse double points, and each double point is provided with an information of being \emph{positive}, \emph{negative} or \emph{virtual} as in Figure~\ref{F:Crossings}. We assume that virtual link diagrams are the same if they are isotopic in~$\R^2$. Positive and negative crossings will also be called \emph{classical crossings}.
\end{defin}

\begin{defin}
\emph{Fused isotopy} is the equivalence relation on the set of virtual link diagrams given by classical Reidemeister moves, virtual Reidemeister moves, the mixed Reidemeister move~\eqref{E:M}, and the forbidden moves \eqref{E:F1} and~\eqref{E:F2}.
\end{defin}

\begin{figure}[htb]
	\centering
		\includegraphics[scale=0.7]{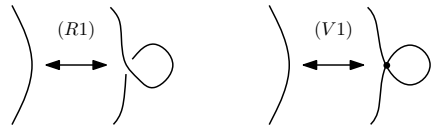} 
	\caption{Reidemeister moves of type I.}
	\label{F:ReidOne}
\end{figure}

\begin{rem}
These moves are the moves pictured in Figure~\ref{F:classical}, \ref{F:virtual}, \ref{F:mixed}, and \ref{F:forbidden}, with the addition of Reidemeister moves of type~I, both classical and virtual, see Figure~\ref{F:ReidOne}. 
\end{rem}

\begin{defin}
A \emph{fused link} is an equivalence class of virtual link diagrams with respect to fused isotopy.
\end{defin}

The classical Alexander Theorem generalizes to virtual braids and links, and it directly implies that every oriented welded (resp.\ fused) link can be represented by a welded (resp.\ unrestricted virtual) braid, whose Alexander closure is isotopic to the original link. Two braiding algorithms are given in~\cite{Kam} and~\cite{KL1}.

Similarly we have a version of Markov Theorem~(\cite{KaL}): before stating it, we recall that the natural map $\UVB\nn \to \UVB{\nn+1}$, that adds one strand on the right of an element of $\UVB\nn$, with the convention of considering braids going from the top to the bottom, is an inclusion. 

\begin{thm}[\cite{KaL}]
Two oriented fused links are isotopic if and only if any two corresponding unrestricted virtual braids differ by moves defined by braid relations in $\UVB\infty$ (braid moves) and a finite sequence of the following moves (extended Markov moves):
\begin{itemize}
\item Virtual and classical conjugation: $\rr\ii \beta \rr\ii \sim \beta \sim \siginv\ii \beta \sig\ii \sim \sig\ii \beta \siginv\ii$;
\item Right virtual and classical stabilization: $\beta \rr\nn \sim \beta \sim \beta \sigma_\nn^{\pm 1}$;
\end{itemize}

where $\UVB\infty = \bigcup_{n=2}^\infty \UVB\nn$, $\beta$ is a braid in $\UVB\nn$, $\sig\ii, \rr\ii$ generators of $\UVB\nn$ and $\sig\nn, \rr\nn \in \UVB{\nn+1}$.
\end{thm}

Here we give an application to fused links of Theorem~\ref{T:uvbn}.

\begin{thm}
\label{T:fused}
Any   fused link  is fused isotopic to the closure of an unrestricted virtual pure braid. 
\end{thm}
 \proof
Let us start remarking that the case of knots is trivial because knots are  fused isotopic to the unknot~(\cite{Nel, Kan}).

Let now $L$ be a fused link with~$\nn>1$ components; then there is an unrestricted virtual braid $\alpha\in \UVB m$ such that $\hat{\alpha}$ is fused isotopic to~$L$.

Let $s_{kl} = \rr{k-1}  \, \rr{k-2} \ldots \rr{l}$ for $l < k$ and $s_{kl} = 1$
in other cases. We define the set 
\[
\Lambda_n = \left\{ \prod\limits_{k=2}^n s_{k,j_k} \vert 1 \leq j_k
\leq k \right\}
\]
which can be seen as the ``virtual part'' of $\UVB\nn$, since it coincides with the set of canonical forms of elements in~$\iota(S_n)$, where $\iota$ is the map from Lemma~\ref{L:lemma1}.

\begin{figure}[htb]
\centering
		\includegraphics[scale=0.7]{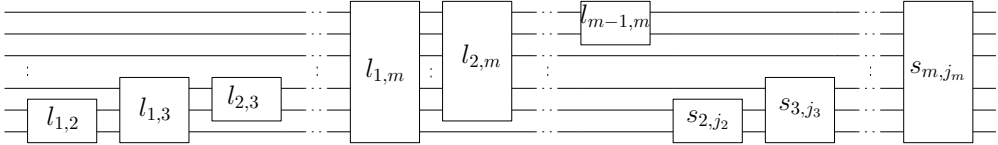} 
	\caption{The braid $\alpha$.}
	\label{F:alpha}
\end{figure}

Then using Theorem~\ref{T:uvbn} we can rewrite $\alpha$ as:
\[
\alpha = l_{1,2} l_{1,3} l_{2,3} \cdots l_{m-1, m} \pi
\]
where $l_{i, j} \in \langle \lam{i,j}, \lam{j,i} \rangle$ and $\pi = s_{2, j_2} \cdots s_{m, j_m} \in \Lambda_n$ (see Figure~\ref{F:alpha}).

Using Lemma~\ref{L:lemma1}, we can do another rewriting:
\[
\alpha = L_{2} s_{2, j_2} L_{3} s_{3, j_3} \cdots L_{m} s_{m, j_m}
\]
where $L_i \in \langle \lam{1,i}, \lam{i, 1} \rangle \times \cdots \times \langle \lam{i-1,i}, \lam{i, i-1} \rangle$.

\begin{figure}[htb]
\centering
		\includegraphics[scale=0.7]{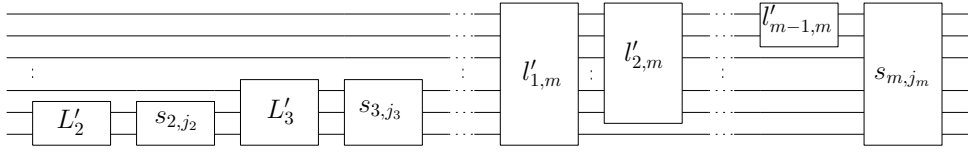} 
	\caption{The rewriting of the braid $\alpha$, with $L^{\prime}_i \in \langle \lambda_{1,i}, \lambda_{i, 1} \rangle \times \cdots \times \langle \lambda_{i-1,i}, \lambda_{i, i-1} \rangle$.}
	\label{F:alpha2}
\end{figure}

Then again we can reorder terms in the~$L_i$s:
\[
\alpha = l^\prime_{1,2} s_{2, j_2} l^\prime_{1,3} l^\prime_{2,3} s_{3, j_3} \cdots l^\prime_{m-1, m} s_{m, j_m}
\]

with $l^\prime_{i, j} \in \langle \lam{i,j}, \lam{j,i} \rangle$, see Figure~\ref{F:alpha2}.

If $s_{i, j_i} =1$ for $i=2, \dots, m$, then $\alpha$ is a pure braid and~$m=n$.

\begin{figure}[htb]
 \centering
		\includegraphics[scale=0.7]{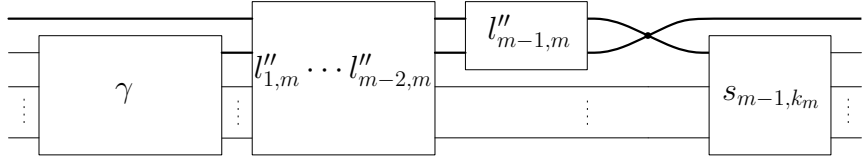} 
	\caption{The form of $\alpha_1$.}
	\label{F:FK}
\end{figure}

Suppose then that there is a $s_{k, j_k} \neq 1$ for some~$k$, and that $s_{i, j_i}=1$ for each~$i > k$. Conjugating $\alpha$ for~$s_{m, 1}^{m-k}$, we obtain a braid $\alpha_1=s_{m, 1}^{k-m} \alpha
s_{m, 1}^{m-k}$ whose closure is fused isotopic to~$L$ where the $k$-th strand of $\alpha$ is the $m$-th strand of~$\alpha_1$. We can rewrite $\alpha_1$ as:
\[
\alpha_1 = \gamma \ l_{1, m}^{\prime\prime} l_{2, m}^{\prime\prime} \cdots l_{m-1, m}^{\prime\prime} \ s_{m, k_m}
\]
 where $\gamma=  l^{\prime\prime} s_{2, j_2}  \cdots l^{\prime\prime}_{m-2, m-1}   s_{m-1, k_{m-1}}$, so it does not involve the $m$-th strand, and $l_{1, m}^{\prime\prime} l_{2, m}^{\prime\prime} \cdots l_{m-1, m}^{\prime\prime}$ is pure. For definition $s_{m, k_m}=\rr{m-1}s_{m-1, {k_m}}$. The $m$-th strand and the other strand involved in this occurrence of $\rr{m-1}$ that we have just isolated, belong to the same component of $L_1=\hat{\alpha_1}$ (see Figure~\ref{F:FK}). Hence also all the crossings in $l_{m-1, m}^{\prime\prime}$ belong to that same component.
 
We virtualize all classical crossings of $l_{m-1, m}^{\prime\prime}$ using Kanenobu's technique~(\cite[Proof of Theorem~1]{Kan}): it consists in deforming the understrand of one classical crossing at a time, considered in the closure of the link, with a sequence of generalized Reidemeister moves, pushing it along the whole component. At the end of the process, there is a new classical crossing instead of the original one, and $2j$ new virtual crossings, where $j$ is the number of crossings the understrand has been pushed through. With generalized Reidemeister moves of braid type, one can change the original classical crossing with a virtual one and remove the new classical crossing with a Reidemeister move of type I. Since our crossings are on the top strand, this Reidemeister move of type I is equivalent to a Markov's classical stabilisation, so we obtain a new link~$L_1^\prime$, fused isotopic to~$L$, associated to a braid $\alpha_1^\prime$ who is identical to $\alpha_1$ except that it has a virtual crossing at the place of the classical crossing considered. This is done for each classical crossing in~$l_{m-1, m}^{\prime\prime}$.

Since $l_{m-1, m}^{\prime\prime}$ has an even total number of generators $\sig{m-1}$ and~$\rr{m-1}$, after virtualizing $l_{m-1, m}^{\prime\prime}\rr{m-1}$ becomes a word composed by an odd number of~$\rr{m-1}$. Applying the relation associated with the virtual Reidemeister move of type~II we obtain a new link~$L_2$, fused isotopic to~$L$, associated to $\alpha_2= \gamma \ l_{1, m}^{\prime\prime} l_{2, m}^{\prime\prime} \cdots  l_{m-2, m}^{\prime\prime} \ \rr{m-1} \ s_{m-1, k_m}$.

Applying once more Lemma~\ref{L:lemma1}, $\alpha_2$ becomes $\gamma \rr{m-1} \  \overline{l_{1, m}}\ \overline{l_{2, m}} \dots  \overline{l_{m-2, m}} \ s_{m-1, k_m}$, where $ \overline{l_{i, m}}$ is a word in~$\langle \lambda_{m-1, \ii} , \lambda_{\ii, m-1}\rangle$.

In $\alpha_2$ there is only one (virtual) crossing on the $m$-th strand, so, using Markov moves (conjugation and virtual stabilisation) we obtain a new braid $\alpha_3$, whose closure is again fused isotopic to $L$ and has $(m-1)$ strands. In other words, the braid $\alpha_3$ is obtained removing from $\alpha_2$ the only virtual crossing on the $m$-th strand, and thanks to Markov theorem its closure is fused isotopic to~$L$.  

If we continue this process, eventually we will get to a braid $\beta$ in $\UVB\nn$ whose closure is fused isotopic to~$L$. At this point, each strand of $\beta$ corresponds to a different component of~$L$, so $\beta$ must be a pure braid. 

 \endproof

The technique used in Theorem~\ref{T:fused} was used, associated with braid decomposition in~$\BB\nn$, by A.~Fish and E.~Keyman to prove the following result about fused links.
\begin{thm}[\cite{FK2}] 
\label{T:FK}
A fused link with only classical crossings $L$ with $n$ components is completely determined by the linking numbers of each pair of components under fused isotopy.
\end{thm}

The proof in \cite{FK2} is quite technical, it involves several computations on generators of the pure braid group and their images in~$\UVP\nn$.
Previous result allows us to give an easier proof: the advantage is that no preliminary lemma on the properties of the pure braid group generators is necessary. 

 \proof[Proof of Theorem~\ref{T:FK}]
We consider a fused link with only classical crossings $L$ with $n$ components: when applying Kanenobu's technique to obtain $\alpha_2$ (see the proof of Theorem~\ref{T:fused}), one gets a braid with only one virtual crossing on the $m$-strand, and removes it, so that the resulting braid $\alpha_3$ only has classical crossings. So, continuing the process, one gets that $L$ is fused isotopic to the closure of an $n$-string unrestricted virtual pure braid~$\beta$ which only has classical crossings.  

Even though $\BB{m}$ and $\PB{m}$ are not subgroups of $\UVB{m}$, since $\hat{\beta}$ has only classical crossings, we can consider $\BB{m}$ and $\PB{m}$'s images in~$\UVB{m}$ and rewrite the pure braid~$\beta$ in terms of $a_{i, j}$ generators, and conclude as Fish and Keyman do, defining a group homomorphism $\delta_{\ii, \jj} : \PB\nn \to \Z$ by 
\[
a_{s, t} \mapsto 
\begin{cases} 1 \quad \mbox{if} \ s=\ii \ \mbox{and} \ t=\jj; \\
0 \quad \mbox{otherwise}
\end{cases}
\]
which is the classical linking number $lk_{i, j}$ of $L$'s $\ii$-th and $\jj$-th components. Any fused link with only classical crossings $L$ with $\nn$ components can be obtained as a closure of a pure braid $\beta = x_2 \cdots x_n$ where each $x_\ii$ can be written in the form $x_\ii=a_{1, \ii}^{\delta_{1, \ii}} \cdots a_{i-1, i}^{\delta_{i-1, i}}$ (Corollary~\ref{C:purecommute}). This shows that $\beta$ only depends on the linking number of the components.

 \endproof

In~\cite[Section~1]{GPV} a \emph{virtual} version of the \emph{linking number} is defined in the following way: to a $2$-component link we associate a couple of integers $(vlk_{1,2}, vlk_{2,1})$ where $vlk_{1,2}$ is the sum of signs of classical crossings where the first component passes over the second one, while $vlk_{2,1}$ is computed by exchanging the components in the definition of $vlk_{1,2}$. Clearly the classical linking number $lk_{1,2}$ is equal to half the sum of $vlk_{1,2}$ and $vlk_{2,1}$.

Using this definition of virtual linking number, we could be tempted to extend Fish and Keyman results, claiming that a fused link $L$ is completely determined by the virtual linking numbers of each pair of components under fused isotopy.

However for the unrestricted case the previous argument cannot be straightforwardly applied: the virtual linking number is able to distinguish $\lam{i,j}$ from $\lam{j,i}$, but it is still an application from $\UVP\nn$ to ${(\Z^2)}^{n(n-1)/2} = \Z^{n(n-1)}$ that counts the exponents (\ie, the number of appearances) of each generator. Since $\UVP\nn$ is not abelian, this is not sufficient to completely determine the braid. 

\begin{rem}
\label{R:linking}
Fish and Keynman in \cite{FK2} suggest that their theorem cannot be extended to links with virtual crossings between different components. They consider the unlink on two components $U_2$ and $L= \hat{\alpha}$, where $\alpha = \sig1 \rr1 \siginv1 \rr1$, they remark that their classical linking number is $0$ but they conjecture that these two links are not fused isotopic. In fact, considering the virtual linking number we can see that $(vlk_{1,2}, vlk_{2,1})(U_2) = (0, 0)$, while~$(vlk_{1,2}, vlk_{2,1})(L) = (-1, 1)$.
\end{rem}


\section{The fused link group}\label{S:fusedgroup}

\subsection{A representation for the unrestricted virtual braid group}

Let us recall that the braid group $\BB\nn$ may be represented as a subgroup of $\Aut(\F\nn)$
by associating to any generator $\sig\ii$, for $\ii=1,2,\ldots,\nno$, of $\BB\nn$ the following
automorphism of $\F\nn$:
\begin{equation}
\label{E:AutoSigma}
\sig\ii : \left\{
\begin{array}{lll}
x_{i} &\longmapsto x_{i} \, x_{i+1} \, x_i^{-1}, &  \\ 
x_{i+1} &\longmapsto
x_{i}, & \\ 
x_{l} &\longmapsto x_{l}, &  l\neq i,i+1.
\end{array} \right.
\end{equation}

Moreover Artin provided (see for instance~\cite[Theorem 5.1]{LH}) a characterization of
braids as automorphisms of free groups: an automorphism
$\beta \in \Aut(\F\nn)$ lies in $\BB\nn$
if and only if $\beta $ satisfies the following conditions:
\begin{enumerate}[label=\roman*)]
\item $\beta(x_i) = a_i \, x_{\pi(i)} \, a_i\inv,~~1\leq i\leq n$\; ;
\item $\beta(x_1x_2 \ldots x_n)=x_1x_2 \ldots x_n$\; ,
\end{enumerate}
where $\pi \in S_\nn$
and $a_i \in \F\nn$.

According to \cite{FRR} we call \emph{group of automorphisms of permutation conjugacy type}, denoted by~$\PC\nn$,
the group of automorphisms satisfying the first condition.
The group~$\PC\nn$ is isomorphic to $\WB\nn$ \cite{FRR}; more precisely to each generator $\sig\ii$ of $\WB\nn$ we associate the previous automorphisms of $\F\nn$ while to each generator~$\rr\ii$, for $i=1,2,\ldots,n-1$, we associate the following automorphism of~$\F\nn$:
\begin{equation}
\label{E:AutoRho}\rr\ii : \left\{
\begin{array}{lll}
x_{i} &\longmapsto   x_{i+1}  &  \\ 
x_{i+1} &\longmapsto
x_{i}, & \\
x_{l} &\longmapsto x_{l}, &  l\neq i,i+1.
\end{array} \right.
\end{equation}

We have thus a faithful representation $\psi \colon \WB\nn \to \Aut(\F\nn)$. 

\begin{rem}
\label{R:manydefinitions}
The group $\PC\nn$ admits also other equivalent definitions in terms of mapping classes and configuration spaces: it
 appears often in the literature with different names and notations, such as group of flying rings~\cite{Bn,BH}, McCool group~\cite{BerP}, motions group~\cite{Go} and loop braid group~\cite{BWC}.
\end{rem}

\begin{rem}
Kamada remarks in \cite{Kam} that the classical braid group $\BB\nn$ embeds in~$\VB\nn$ through the canonical epimorpism~$\VB\nn \to \WB\nn$. It can be seen via an argument in \cite{FRR} that $\BB\nn$ is isomorphic to the subgroup of $\VB\nn$ generated by~$\{ \sig1, \dots, \sig\nn \}$.
\end{rem}
 

\begin{rem} 
\label{R:F2}
As a consequence of the isomorphism between $\WB\nn$ and~$\PC\nn$, we can show that relation \eqref{E:F2} does not hold in~$\WB\nn$. In fact applying $\rr{\ii+1} \sig\ii \sig{\ii+1}$ one gets
\begin{align*}
\rr{\ii+1} \sig\ii \sig{\ii+1}&:
\left\{
\begin{array}{ll}
x_\ii  \longmapsto x_\ii  \longmapsto x_\ii x_{\ii+1} x_\ii\inv \longmapsto x_\ii x_{\ii+1} x_{\ii+2} x_{\ii+1}\inv x_\ii\inv, \\
x_{\ii+1}   \longmapsto x_{\ii+2}  \longmapsto x_{\ii+2}  \longmapsto x_{\ii+1}, \\
x_{\ii+2}  \longmapsto x_{\ii+1}  \longmapsto x_\ii  \longmapsto x_\ii  ,
\end{array} \right.
\intertext{
while applying $\sig\ii \sig{\ii+1} \rr{\ii}$ one gets}
\sig\ii \sig{\ii+1} \rr{\ii}&:
\left\{
\begin{array}{ll}
x_\ii  \longmapsto x_\ii x_{\ii+1} x_\ii\inv   \longmapsto  x_\ii x_{\ii+1} x_{\ii+2} x_{\ii+1}\inv x_\ii\inv \longmapsto x_{\ii+1} x_\ii x_{\ii+2} x_\ii\inv x_{\ii+1}\inv, \\
x_{\ii+1}   \longmapsto x_\ii  \longmapsto x_\ii \longmapsto x_{\ii+1}, \\
x_{\ii+2}  \longmapsto x_{\ii+2}  \longmapsto x_{\ii+1} \longmapsto x_\ii .
\end{array} \right. 
\end{align*}
Since $ x_\ii x_{\ii+1} x_{\ii+2} x_{\ii+1}\inv x_\ii\inv \neq x_{\ii+1} x_\ii x_{\ii+2} x_\ii\inv x_{\ii+1}\inv$ in $\F\nn$ we deduce that relation \eqref{E:F2} does not hold in~$\WB\nn$.
\end{rem}

Our aim is to find a representation for unrestricted virtual braids as automorphisms of a group~$G$. Since the map $\psi \colon \WB\nn \to \Aut(\F\nn)$ does not factor through the quotient $\UVB\nn$ (Remark~\ref{R:F2})
we need to find a representation in the group of automorphisms of a quotient of $\F\nn$ in which relation \eqref{E:F2} is preserved.

\begin{rem}
In \cite{KMRW} the authors look for representations of the braid group $\BB\nn$ that can be extended to the loop braid group $\WB\nn$ but do not factor through $\UVB\nn$, which is its quotient by relations of type \eqref{E:F2}, while we look for a representation that does factor. 
\end{rem}

Let $\F\nn =  \gamma_1 \F\nn \supseteq \gamma_2 \F\nn \supseteq \cdots $ be the lower central series of $\F\nn$, the free group of rank~$\nn$, where~$\gamma_{i+1}\F\nn=[\F\nn, \gamma_i \F\nn]$.
Let us consider its third term, $\gamma_3 \F\nn=\big[\F\nn, [\F\nn, \F\nn]\big]$; the free $2$-step nilpotent group $N_\nn$ of rank $\nn$ is defined to be the quotient~$\faktor{\F\nn}{\gamma_3\F\nn}$.

There is an epimorphism from $\F\nn$ to $N_\nn$ that induces an epimorphism from $\Aut(\F\nn)$ to $\Aut(N_\nn)$~(see \cite{And}).
Then, let $\phi \colon \UVB\nn \to \Aut(N_\nn)$ be the composition of $\varphi \colon \UVB\nn \to \Aut(\F\nn)$ and~$\Aut(\F\nn) \to \Aut(N_\nn)$.

\begin{prop}
\label{P:representation}
The map $\phi \colon \UVB\nn \to \Aut(N_\nn)$ is a representation for $\UVB\nn$. 
\end{prop}

 \proof
We use the convention $[x, y] = x\inv y\inv x y$. In $N_\nn$ we have that 
$ \big[[x_\ii, x_{\ii+1}], x_{\ii+2}\big] = 1$, for $\ii = 1, \dots, \nn-2$, meaning that 
$ x_\ii x_{\ii+1} x_{\ii+2} x_{\ii+1}\inv x_\ii\inv = x_{\ii+1} x_\ii x_{\ii+2} x_\ii\inv x_{\ii+1}\inv$, \ie, relation \eqref{E:F2} is preserved.

 \endproof

\begin{prop}
\label{P:Immagine}
The image of the representation $\phi \colon \UVP\nn \to \Aut(N_\nn)$ is a free abelian group of rank~$\nn(\nno)$.
\end{prop}
 \proof
From Theorem~\ref{T:uvbn} we have that the only generators that do not commute in $\UVP\nn$ are $\lam{\ii, \jj}$ and $\lam{\jj ,\ii}$ with $1 \leq \ii \neq \jj \leq \nn$.

Recalling the expressions of $\lam{\ii ,\jj}$ and $\lam{\jj ,\ii}$ in terms of generators $\sig\ii$ and $\rr\ii$, we see that the automorphisms associated to $\lam{\ii, \jj}$ and $\lam{\jj ,\ii}$ are
\begin{align*}
\phi(\lam{\ii, \jj})&: 
\left\{
\begin{array}{ll} 
x_\ii  \longmapsto x_\jj\inv x_\ii x_\jj = x_\ii [x_\ii, x_\jj] \\
x_k \longmapsto x_k, \ \mbox{for} \ k \neq \ii;
\end{array} \right. \\
\phi(\lam{\jj ,\ii})&:
\left\{
\begin{array}{ll} 
x_\jj \longmapsto  x_\ii\inv x_\jj x_\ii = x_\jj[x_\jj, x_\ii] = x_\jj[x_\ii, x_\jj]\inv;\\
x_k \longmapsto x_k \ \mbox{for} \ k \neq \ii.
\end{array} \right.
\intertext{It is then easy to check that the automorphisms associated to $\lam{\ii, \jj}\lam{\jj, \ii}$ and to $ \lam{\jj, \ii} \lam{\ii, \jj}$ coincide:}
\phi(\lam{\ii, \jj} \lam{\jj, \ii} )=\phi( \lam{\jj, \ii} \lam{\ii, \jj})&: 
\left\{
\begin{array}{ll} 
x_\ii  \longmapsto x_\ii [x_\ii, x_\jj]\\
x_\jj \longmapsto x_\jj[x_\ii, x_\jj]^{-1}.
\end{array} \right. \\
\intertext{To see that in $\phi(\UVP\nn)$ there is no torsion, let us consider a generic element $w$ of~$\UVP\nn$. It will have the form $w= l_{1,2} l_{1,3} \cdots l_{n-1, n}$ where $l_{i,j}$ is a product of generators $\lam{i,j}$ and~$\lam{j,i}$.
Generalizing the calculation done above, we have that}
\phi(l_{1,2} l_{1,3} \cdots l_{n, n-1})&=\phi(\lamm{1,2}{\varepsilon_{1,2}} \lamm{2,1}{\varepsilon_{2,1}} \cdots \lamm{n-1, n}{\varepsilon_{n-1, n}} \lamm{n, n-1}{\varepsilon_{n, n-1}}),
\intertext{where $\varepsilon_{i,j}$ is the total number of appearances of $\lam{i,j}$ in~$l_{i,j}$.
With another easy calculation (check out also Remark~\ref{R:calculation}) we have that:}
\phi(\lamm{1,2}{\varepsilon_{1,2}} \lamm{2,1}{\varepsilon_{2,1}} \cdots \lamm{n-1, n}{\varepsilon_{n-1, n}} \lamm{n, n-1}{\varepsilon_{n, n-1}})&: 
\left\{
\begin{array}{llll} 
x_1  \longmapsto x_1 [x_1, x_2]^{\varepsilon_{12}} [x_1, x_3]^{\varepsilon_{13}} \cdots [x_1, x_n]^{\varepsilon_{1n}}  \\
x_2  \longmapsto x_2 [x_2, x_1]^{\varepsilon_{21}} [x_2, x_3]^{\varepsilon_{23}} \cdots [x_2, x_n]^{\varepsilon_{2n}}  \\
\vdots \\
x_n  \longmapsto x_n [x_n, x_1]^{\varepsilon_{n1}} [x_n, x_2]^{\varepsilon_{n2}} \cdots [x_n, x_{n-1}]^{\varepsilon_{n, {n-1}}} \\
\end{array} \right. \\
\end{align*}
So the condition for $\phi(w)$ to be $1$ is that all exponents are equal to~$0$, hence~$w=1$.

 \endproof

\begin{rem}
As a consequence of the previous calculation the homomorphism $\phi$ coincides on $\UVP\nn$ with the abelianization map.
\end{rem}

As a consequence of Proposition~\ref{P:Immagine}, the representation $\phi$ is not faithful.
However, according to the previous characterization of $WB_n$ as subgroup of $\Aut(\F\nn)$ 
it is natural to ask if we can give a characterization of automorphisms of $\Aut(N_\nn)$ that belong to~$\phi(\UVB\nn)$.

\begin{prop}\label{P:caractaut} 
Let $\beta$ be an element of~$\Aut(N_\nn)$, then $\beta \in  \phi(\UVB\nn)$ if and only if $\beta$ satisfies the condition $\beta(x_i) = a_i^{-1}  x_{\pi(i)}  a_i$ with~$1\leq \ii\leq n$, where $\pi \in S_\nn$ and~$a_i \in N_\nn$.
\end{prop}

 \proof
Let us denote with $UVB(N_n)$ the subgroup of $\Aut(N_n)$ such that any element $\beta \in UVB(N_n)$ has the form $\beta(x_i)=  g_i^{-1}  x_{\pi(i)} g_i$, denoted by~$x_{\pi(i)}^{g_i}$, with $1\leq \ii\leq n$, where $\pi \in S_\nn$ and~$g_i \in N_\nn$. We need to prove that $\phi \colon \UVB\nn \to UVB(N_n)$ is an epimorphism. 
Let $\beta$ be an element of~$UVB(N_n)$. Since $S_n$ is both isomorphic to the subgroup of $\UVB\nn$ generated by the $\rr\ii$ generators, and to the subgroup of $UVB(N_n)$ generated by the permutation automorphisms, we can assume that for $\beta$ the permutation is trivial, \ie,~$\beta(x_i) = x_\ii^{g_i}$.
We define $\varepsilon_{i,j}$ to be $\phi(\lam{i,j})$ as in Proposition~\ref{P:Immagine}, and we prove that $\beta$ is a product of such automorphisms. 
We recall that $x^{ y z}=x^{z y}$ for any~$x,y, z \in N_n$, therefore: 

\[
\beta(x_i) = x_\ii^{x_1^{a_{i,1}} \cdots x_n^{a_{i,n}}}
\] 
where~$a_{i,i}=0$.

In particular we can assume that
\[
\beta(x_1) = x_1^{x_2^{a_{1,2}} \cdots \, x_n^{a_{1,n}}}.
\]

We define a new automorphism $\beta_1$ multiplying $\beta$ by $\varepsilon_{1,2}^{-a_{1,2}} \cdots \varepsilon_{1,n}^{-a_{1,n}}$. 
We have that $\beta_1(x_1)=x_1$, and $\beta_1(x_j)=\beta(x_j)$ for $j \not= 1$. Then again we define a new automorphism $\beta_2= \beta_1 \ \varepsilon_{1,2}^{-a_{2,1}} \varepsilon_{2,3}^{-a_{2,3}} \cdots \varepsilon_{2, n}^{-a_{2,n}}$ that fixes $x_1$ and $x_2$.

Carrying on in this way for $\nn$ steps we get to an automorphism 
\[
\beta_n = \beta_{n-1} \ \varepsilon_{n, 1}^{-a_{n,1}} \cdots  \varepsilon_{n, n-1}^{-a_{n,n-1}} = \beta  \ \prod_{j=1}^{n} \varepsilon_{n, j}^{-a_{1,j}} \prod_{j=1}^{n} \varepsilon_{n-1, j}^{-a_{2,j}} \cdots \prod_{j=1}^{n} \varepsilon_{1, j}^{-a_{n,j}} 
\] 
setting $\varepsilon_{i,i}=1$. The automorphism $\beta_n$ is the identity automorphism: then $\beta$ is a product of $\varepsilon_{i,j}$ automorphisms, hence it has a pre-image in $\UVB\nn$.

 \endproof


\subsection{The fused link group}

Let $L$ be a fused link. Then there exists an unrestricted virtual braid $\beta$ such that its closure $\hat{\beta}$ is equivalent to~$L$. 
\begin{defin}
The fused link group $G(L)$ is the group given by the presentation
\[
\bigg\langle \  x_1, \dots, x_\nn \ \bigg\vert \
\begin{matrix}
\phi(\beta)(x_\ii)=x_\ii &\mbox{ for } \ii \in \{1, \dots, \nn \}, \\
\big[x_i, [x_k,  x_l] \big]=1 &\mbox{ for } i, k, l \mbox{ not necessarily distinct} \end{matrix} \
\bigg\rangle
\]
where $\phi \colon \UVB\nn \to \Aut(N_\nn)$ is the map from Proposition~\ref{P:representation}. 
\end{defin}

\begin{prop}
\label{P:Group}
The fused link group is invariant under fused isotopy.
\end{prop}
 \proof
  According to~\cite{KaL} two unrestricted virtual braids have fused isotopic closures if and only if they are related by \emph{braid moves} and \emph{extended Markov moves}. We should check that under these moves the fused link group $G(L)$ of a fused link $L$ does not change. This is the case. However a quicker strategy to verify the invariance of this group is to remark that it is a projection of the \emph{welded link group} defined in~\cite[Section~5]{BB2}. This last one being an invariant for welded links, we only have to do the verification for the second forbidden braid move, coming from relation~\eqref{E:F2}. This invariance is guaranteed by the fact that $\phi$ preserves relation~\eqref{E:F2} as seen in Proposition~\ref{P:representation}.

 \endproof

\begin{figure}[htb]
	\centering
		\includegraphics[scale=0.7]{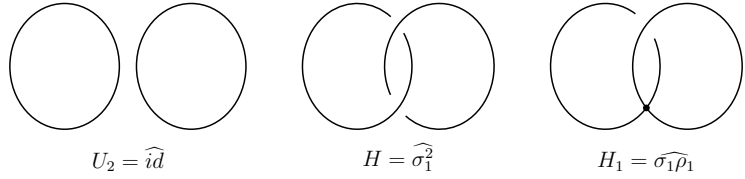} 
	\caption{The fused link group distinguishes the unlink $U_2$ from the Hopf link $H$, but does not distinguish the Hopf link with two classical crossings $H$ from the one with a classical and a virtual crossing~$H_1$. In fact: $G(U_2) = N_2$, while $G(H) = G(H_1) = \Z^2$. We remark however that $H$ and $H_1$ are distinguished by the virtual linking number.}
	\label{F:Hopf}
\end{figure}

\begin{rem}
\label{R:calculation}
Let us recall that, according Theorem \ref{T:fused},  a fused link $L$ admits as a representative the closure of an element of $UVP_n$, say $\beta_L$, and
following  the proof of Proposition \ref{P:caractaut}, we can deduce that 
\[
\phi(\beta(x_i)) = x_\ii^{x_1^{a_{i,1}} \cdots x_n^{a_{i,n}}}
\] 
where~$a_{i,i}=0$ and   $a_{i,j}=vlk_{i,j}$ for $i\not=j$;  Since virtual linking numbers are fused invariants, we get another easy proof 
of Proposition \ref{P:Group}. However, it means also that the knot group is determined by  virtual linking numbers; actually, as shown in Figure~\ref{F:Hopf},
is weaker.  The  relation between   virtual linking numbers and the knot group 
can be nicely described in the case $n=2$ as follows. 
Let us consider $\lamm{1,2}\alpha \lamm{2,1}\beta$ and~$\lamm{1,2}\gamma$, where~$\gamma$ is the greatest common divisor of $\alpha$ and~$\beta$
and therefore of  $vlk_{1,2}$ and $vlk_{2,1}$. 
The automorphisms associated to them are 

\begin{align*}
\phi( \lamm{1, 2}\alpha \lamm{2, 1}\beta)&: 
\left\{
\begin{array}{ll} 
x_1  \longmapsto x_1 \big[x_1, x_2 [x_1, x_2]^{-\beta}\big]^\alpha = x_1 \big[x_1, [x_1, x_2]^{-\beta}\big]^\alpha [x_1, x_2]^{\alpha } =x_1 [x_1, x_2]^\alpha\\
x_2 \longmapsto x_2[x_1, x_2]^{-\beta} ;
\end{array} \right. \\
\phi(\lamm{1,2}\gamma)
&:  \left\{
\begin{array}{ll} 
x_1  \longmapsto x_1 [x_1, x_2]^\gamma\\
x_2 \longmapsto x_2 .
\end{array} \right.
\end{align*}


Then
\[G(\lamm{1, 2}\alpha \lamm{2, 1}\beta) = G(\lamm{1,2}\gamma) = \Gr{x_1, x_2 \mid [x_1, x_2]^{\gamma}=1, \ \big[x_i, [x_k,  x_l] \big]=1 \mbox{ for } i, k, l \in \{1, 2\}}= \\
\]
\[
= \Gr{x_1, x_2, t \mid [x_1, x_2]=t, t^{\gamma}=1,  t \; \mbox{central} }\]
This latter group presentation allows to distinguish these groups for different $\gamma \in \N$ (since~$\gamma$  is the order of the central element $t$ of these Heisenberg-like groups,
setting that $\gamma=0$ means that $t$ has infinite order); in particular  we can set $G(\lamm{1,2}\gamma):=G_\gamma$. For instance  
 the two links considered in \cite{FK1}, $L= \widehat{\sig1 \rr1 \siginv1 \rr1}$ and $U_2$, have corresponding groups   $G_1=\Z^2$ and $G_0=N_2$ and therefore are distinguished
 by  $G_\gamma$, while, as we saw above,  have the same classical linking number.
\end{rem}


\section{Other quotients} \label{S:quotients}

Several other quotients of virtual braid groups have been studied in the literature: we end this paper with a short survey
on them, giving the structure of the corresponding pure subgroups and some results on their linearity.

\subsection{Flat virtual braids.} \label{S:Flat}

The study of flat virtual knots and links was initiated by Kauffman \cite{Ka} and their braided counterpart was introduced in~\cite{Ka1}. The category of flat virtual knots is identical to the structure of what are called virtual strings by V.~Turaev in~\cite{T} (remark that every virtual string is the closure of a flat virtual braid).

The flat virtual braid grous $FVB_\nn$ was introduced in \cite{Ka1} as a
quotient of~$\VB\nn$ adding relations
\begin{equation}
\label{E:quadratic}
\sigg\ii2 = 1\mbox{,} \quad \mbox{for } \ii = 1, \dots, \nno.
\end{equation}

It is evident that $FVB_\nn$ is a quotient of the free product~$S_\nn \ast S_\nn$.

Let us consider the natural projection map $f \colon \VB\nn \to FVB_\nn$, and set $f(\rr\ii):=\rr\ii$ and~$f(\sig\ii):=s_\ii$ for~$\ii = 1, \dots, \nno$.

In addition to relations coming from the two copies of~$S_\nn$, in $FVB_\nn$ we have mixed relations 
\begin{align}
\label{E:rel1} s_\ii \rr\jj = \rr\jj s_\ii, \quad  &\mbox{for } \vert \ii - \jj \vert \geq 2, \\
\label{E:rel2} \rr\ii \rr{\ii+1} s_\ii = s_{\ii+1} \rr\ii \rr{\ii+1}, \quad &\mbox{for }\ii = 1, \dots, n-2.
\end{align}

We call \emph{flat virtual pure braid group} $FVP_\nn$ the kernel of the map $FVB_n \to S_n$ defined by $s_i, \rr\ii \mapsto (i, i+1)$ for~$\ii = 1, \dots, \nno$. With respect to the map~$f \colon VB_\nn \to FVB_\nn$, we have that~$f(VP_\nn)= FVP_\nn$.

\begin{prop}
\label{P:VirtualFlat}

Let $VP_\nn^+$ be the (abstract) presented group
\[
\bigg\langle \ \{\lam{i,j} \mid 1 \leq \ii < \jj \leq \nn\} \  \bigg\vert \
\begin{matrix}
\lam{k,l} = \lam{k,l} \lam{i,j},\\
\lam{k,i}  (\lam{k,j} \lam{i,j}) = (\lam{i,j} \lam{k,j})   \lam{k,i} \end{matrix} \
\bigg\rangle.
\]
Then $\VP\nn^+$ coincides with the subgroup of $VP_\nn$ generated by the set $\{\lam{\ii,\jj} \mid 1 \leq \ii < \jj \leq \nn\}$ and is isomorphic to~$FVP_\nn$.
\end{prop}

 \proof

First let us recall that $VP_\nn$ is generated by elements $\lam{i,j}$ defined in Eq.~\eqref{E:lambda} and~\eqref{E:lambda2}, and has the following complete set of relations:
\begin{align}
\tag{RS1} \qquad  \lam{i,j}   \lam{k,l} &= \lam{k,l}   \lam{i,j} ,\\
\tag{RS2} \qquad \lam{k,i}  (\lam{k,j} \lam{i,j}) &= (\lam{i,j} \lam{k,j})   \lam{k,i}.
\end{align}

Now define  the map $\iota \colon \VP\nn^+ \to \VP\nn$ sending $\lam{i,j}$ to $\lam{i,j}$ and the map
 $\theta \colon \VP\nn \to \VP\nn^+$ sending $\lambda_{\ii ,\jj}$ to $\lambda_{\ii,\jj}$ if $\ii < \jj$ or to  $\laminv{\jj,\ii}$ whenever $\ii > \jj$. Both $\iota$ and $\theta$ are well defined
 homorphisms and $\theta \circ \iota= Id_{ \VP\nn^+ }$ so $\iota$ is injective.

Setting $f(\lam{i,j})=\mu_{i,j}$ and proceeding with similar arguments as in  Theorem \ref{T:theorem1} one can easily prove that  $FVP_\nn$ admits the presentation: 
\[
\Bigg\langle \ \{\mu_{i,j} \mid 1 \leq \ii \neq \jj \leq \nn\} \  \Bigg\vert \
\begin{matrix}
\mu_{i,j} \mu_{k,l} = \mu_{k,l} \mu_{i,j},\\
\mu_{k,i}  (\mu_{k,j} \mu_{i,j}) = (\mu_{i,j} \mu_{k,j})   \mu_{k,i} , \\
\mu_{i, j} \mu_{j, i} =1 \quad \mbox{for } 1 \leq \ii \leq \nno
\end{matrix} \
\Bigg\rangle.
\]

We can proceed as before and to consider the abstract group $FVP_\nn^+$ given by following presentation:
\[
\bigg\langle \ \{\mu_{i,j} \mid 1 \leq \underline{\ii < \jj} \leq \nn\} \  \bigg\vert \
\begin{matrix}
\mu_{i,j} \mu_{k,l} = \mu_{k,l} \mu_{i,j}\\
\mu_{k,i}  (\mu_{k,j} \mu_{i,j}) = (\mu_{i,j} \mu_{k,j})   \mu_{k,i} \\
\end{matrix} \
\bigg\rangle.
\]
We can therefore consider map $\iota' \colon FVP_\nn^+ \to FVP_\nn$ sending $\mu_{i,j}$ to $\mu_{i,j}$ and the map
 $\theta' \colon FVP_\nn \to FVP_\nn^+$ sending $\mu_{\ii ,\jj}$ to $\mu_{\ii,\jj}$ if $\ii < \jj$ or to  $\mu_{\jj,\ii}^{-1}$ whenever $\ii > \jj$. Both $\iota'$ and $\theta'$ are well defined
 homeomorphisms and $\theta' \circ \iota'= Id_{ FVP_\nn^+ }$ and $\iota' \circ \theta'= Id_{ FVP_\nn }$.  
Then  $FVP_\nn^+$ is a group presentation for $FVP_\nn$ and  the isomorphism of the statement is obviously obtained  sending $\mu_{i,j}$ to $\lam{i,j}$.

 \endproof

\begin{rem}
For $\nn=3$, the group 
\[
FVP_3 = \Gr{
\lam{1,2}, \lam{1,3}, \lam{2,3} \ \mid \ \laminv{1,2} (\lam{2,3} \lam{1,3} ) \lam{1,2} =
\lam{1,3} \lam{2,3}
}
\]
is the HNN-extension of the free group $\Gr{ \lam{1,3}, \lam{2,3}}$ of rank $2$ with stable element
$\lam{1,2}$ and with associated subgroups $A = \langle \lam{2,3} \lam{1,3} \rangle$ and
$B = \langle \lam{1,3} \lam{2,3} \rangle$, which are isomorphic to the infinite cyclic group.
Moreover, the group $FVP_3$ is isomorphic to the free product~$ \mathbb{Z}^2 \ast \mathbb{Z}$.
The first claim follows from the previous theorem.
The second one  follows from the observation that setting $a = \lambda_{23} \lambda_{13}$, $b = \lambda_{23}$
we obtain the following new presentation:
\[
FVP_3 = \langle \lambda_{12}, a, b~||~\lambda_{12}^{-1} a \lambda_{12} =
b^{-1} a b \rangle;
\]
if we denote $c = b \lambda_{12}^{-1}$ and exclude $\lambda_{12}$ from the set of generators we get
\[
FVP_3 = \langle  a, b, c~||~[a, c] =
1 \rangle = \langle a, c~|~[a, c] = 1 \rangle * \langle b \rangle.
\]
\end{rem}

Let us recall that there is another remarkable surjection
of the virtual braid group $VB_n$ onto the symmetric group~$S_n$,
which sends $\sig\ii$ into $1$ and $\rr\ii$ into~$(i, i+1)$: the kernel 
of this map is denoted by $H_n$ in~\cite{BB}.
In the same way 
 we can define the group $FH_n$ as the kernel of the homomorphism
$\mu \colon FVB_n \to S_n$, which is defined as follows:
\[
\mu(s_i)=1, \; \mu(\rr\ii)=(i, i+1), \;i=1,2,\dots, n-1.
\]
Now let us define, for $\ii = 1, \dots, \nno$:
\begin{equation}
\begin{split}
 y_{i,i+1}&=s_i, \\ 
 y_{i+1,i}&=\rr\ii s_i \rr\ii.
\end{split}
\end{equation}
For $1 \leq \ii < \jj-1 \leq \nno$:
\begin{equation}
\begin{split}
y_{i,j}&=\rr{j-1} \cdots \rr{i+1} s_i \rr{i+1} \cdots \rr{j-1}, \\
y_{j,i}&=\rr{j-1} \cdots \rr{i+1} \rr\ii s_i \rr\ii \rr{i+1} \cdots \rr{j-1}.
\end{split}
\end{equation}
It is not difficult to prove that these elements belong to $FH_n$ and that:

\begin{thm} \label{T:flatH}
The group $FH_n$ admits a presentation with generators $y_{k,\, l},$ for $1 \leq k \neq l \leq
n$,
and defining relations:
\begin{align}
y_{k,l}^2 &= 1,\\
y_{i,j} \,  y_{k,\, l} &= y_{k,\, l}  \, y_{i,j} \Leftrightarrow (y_{i,j} \,  y_{k,\, l})^2 = 1,\\
y_{i,k} \,  y_{k,j} \,  y_{i,k} &=  y_{k,j} \,  y_{i,k} \, y_{k,j} \Leftrightarrow (y_{i,k} \,  y_{k,j} )^3 = 1,
\end{align}
where distinct letters stand for distinct indices.
\end{thm}

\proof
We can use  Reidemeister-Schreier method and check the above set of relations is complete
or we can consider  a commutative diagram similar to the one of proof of  Theorem~\ref{T:theorem1}:
\[
\begin{CD}
@. 1 @. 1 \\
@. @VVV  @VVV \\
@. \ker f_{\vert H_n} @>>> \ker f \\
@. @VVV  @VVV @. \\
1 @>>> H_n @>>>  \VB\nn @>>> S_\nn @>>>1 \\
   @.  @VV{ f_{\vert H_n}}V         @VV{f}V        @|  \\
1 @>>> FH_n @>>>  FVB_n @>>> S_\nn @>>>1 \\
@. @VVV  @VVV \\
@. 1 @. 1 \\
\end{CD}
\] 
Recall also that, according to~\cite{BB},
 the group $H_n$ is generated by following elements:
\begin{equation}
\begin{split}
& h_{i,i+1}=\sig\ii, \\ 
& h_{i+1,i}=\rr\ii \sig\ii \rr\ii ,
\end{split}
\end{equation}

and for $1 \leq \ii < \jj-1 \leq \nno$:
\begin{equation}
\begin{split}
& h_{i,j}=\rr{j-1} \cdots \rr{i+1} \sig\ii \rr{i+1} \cdots \rr{j-1}, \\
& h_{j,i}=\rr{j-1} \cdots \rr{i+1} \rr\ii \sig\ii \rr\ii \rr{i+1} \cdots \rr{j-1},
\end{split}
\end{equation}

with defining relations:
\begin{align} 
h_{i,j} \,  h_{k,\, l} &= h_{k,\, l}  \, h_{i,j} ,\\
h_{i,k} \,  h_{k,j} \,  h_{i,k} &=  h_{k,j} \,  h_{i,k} \, h_{k,j},
\end{align}
where distinct letters stand for distinct indices.
Now, remarking first that $f(h_{i,j})=y_{i,j}$,  $\ker f = \ker f_{\vert H_n}= \ll \sigg\ii{2} \vert \; i = 1, 2, \ldots, n-1\gg$
and 
$\sigg\ii{2}=h_{i,i+1}^2$, one can also verify that
 $\ker  f_{\vert H_n}$ is generated by elements of type  $g h_{k,l}^2 g^{-1}$, for $1 \le k\not=l \le n$
 and $g\in H_n$ (details are left to the reader, but arguments are the same as in Theorem ~\ref{T:theorem1}).
 Therefore we have the expected complete set of relations for~$FH_n$.
 
 \endproof

\begin{cor}
The group $FVB_n$ is linear.
\end{cor}

 \proof
From the decomposition $VB_n = H_n \rtimes S_n$ we have that~$FVB_n = FH_n \rtimes S_n$,
where $FH_n$ is a finitely generated Coxeter group.
The statement therefore follows from 
the fact that all finitely generated Coxeter groups are linear and that finite extensions of linear groups are also linear.

 \endproof

\subsection{Flat welded braids}
In a similar way we can define the flat welded braid group $FWB_\nn$ as the 
quotient of~$\WB\nn$ adding relations
\begin{equation}
\label{E2:quadratic}
\sigg\ii2 = 1\mbox{,} \quad \mbox{for } \ii = 1, \dots, \nno.
\end{equation}

Let us consider the natural projection map $g \colon \VB\nn \to FVB_\nn$, and set $g(\rr\ii)=\rr\ii$ and~$g(\sig\ii)=s_\ii$ for~$\ii = 1, \dots, \nno$.

In~$FWB_n$, in addition to relations \eqref{E:rel1} and~\eqref{E:rel2}, we also have relations coming from relations of type \eqref{E:F1}, \ie,
\begin{equation}\label{E:rel3}
s_{\ii+1} s_{\ii} \rr{\ii+1} = \rr\ii s_{\ii+1} s_{\ii},\quad \mbox{for }\ii = 1, \dots, \nno.
\end{equation}

In $FWB_\nn$ relations \eqref{E2:quadratic} and \eqref{E:rel3} imply that also relations of type \eqref{E:F2} hold, since from $\rr\ii s_{\ii+1} s_\ii = s_{\ii+1} s_\ii \rr{\ii+1}$ one gets $s_\ii s_{\ii+1} \rr\ii = \rr{\ii+1} s_\ii s_{\ii+1}$.

Adapting Theorem~\ref{T:theorem1} one can easily verify that $FWP_\nn$ is isomorphic to~$\Z^{\nn(\nno)/2}$.
As a straightforward consequence of Theorem \ref{T:uvbn}, we can describe the structure of~$FWB_\nn$.

\begin{prop}
\label{P:fwbn}
Let $\Z^{n(n-1)/2}$ be the free abelian group of rank~$\nn(\nno)/2$. Let us denote by $x_{\ii,\jj}$ for $1 \le \ii \neq \jj \le \nn$
a possible set of generators.
The group $FWB_\nn$ is isomorphic to~$\Z^{n(n-1)/2} \rtimes S_\nn$, where $S_\nn$ acts by permutation on the indices of generators
of $\Z^{n(n-1)/2}$ (setting $x_{\jj,\ii}:=x_{\ii,\jj}\inv$ for~$1 \le \ii < \jj \le \nn$).
\end{prop}

 \proof
Let us recall how elements $\lam{i,j}$ in $\UVB\nn$ were defined.

For $\ii = 1, \dots, \nno$:
\begin{align*}
&\lam{i,i+1} = \rr\ii \siginv\ii,\\
&\lam{i+1,i} = \rr\ii \lam{i,i+1} \rr\ii = \siginv\ii \rr\ii.
\end{align*}

For $1 \leq \ii < \jj-1 \leq \nno$:
\begin{align*}
&\lam{i,j} = \rr{j-1} \rr{j-2} \ldots \rr{i+1} \lam{i,i+1} \rr{i+1} \ldots \rr{j-2} \rr{j-1}, \\
&\lam{j,i} = \rr{j-1} \rr{j-2} \ldots \rr{i+1} \lam{i+1,i}  \rr{i+1} \ldots \rr{j-2} \rr{j-1}.
\end{align*}

Relations \eqref{E2:quadratic} are therefore equivalent to relations~$\lam{i,j} \lam{j,i} =1$.
 Adding these relations and following verbatim the proof of Theorem \ref{T:theorem1} we get the statement.
 
  \endproof
  
\medskip


\subsection{Virtual Gauss braids.}

 From the notion of flat virtual knot we can get the notion of Gauss virtual knot or simply Gauss knot. Turaev \cite{T1}
introduced these knots under the name of ``homotopy classes of Gauss words'', while Manturov \cite{Man1} used the name ``free knots''.

The ``braided'' analogue of Gauss knots, called \emph{free virtual braid group on $n$ strands}, was  introduced in~\cite{Man}. 
From now on we will be calling it \emph{virtual Gauss braid group} and will denote it by~$GVB_n$.

The group of \emph{virtual Gauss braids} $GVB_n$ is the quotient of $FVB_n$ by relations
\[
s_i \rr\ii = \rr\ii s_i \mbox{,}  \quad \mbox{for} \quad i = 1, \dots, n-1.
\]

Note also that the virtual Gauss braid group is a natural quotient of the \emph{twisted virtual braid group}, studied for instance in~\cite{L}.

Once again we can consider the homomorphism from $GVB_\nn$ to $S_\nn$ that sends each generator $s_\ii$ and $\rr\ii$ in~$\rr\ii$. The \emph{virtual Gauss pure braid group} $GVP_\nn$ is defined to be the kernel of this map. Since this map admits a natural section $GVB_n$ is isomorphic to~$GVP_n \rtimes S_n$.

Adapting the proof of Theorem~\ref{T:theorem1}, we get the following.

\begin{prop}
The group $GVP_n$ admits a presentation with generators~$\lam{k,l}$ for $1 \leq k < l \leq
n$ and the defining relations of $FVP_n$ plus relations
\[
\lamm{i,j}2 = 1 \mbox{,} \quad \mbox{for} \quad  1 \leq i < j \leq n.
\]
\end{prop}

Moreover as in the case of $FVB_n$ also in the case of $GVB_n$
we can consider the map~$\mu: GVB_n \to S_n$, defined as follows:
\[
\mu(s_i)=1 \mbox{, } \mu(\rr\ii)=\rr\ii \mbox{,} \quad \mbox{for} \quad i=1,2,\dots, n-1.
\]

Let $GH_n$ be the kernel of the map $\mu: GVB_n \to S_n$, and $y_{k, l},$ the elements defined in subsection~\ref{S:Flat}: we can prove the following result.

\begin{thm} \label{T:gaussH}
The group $GH_n$ admits a presentation with generators $y_{k, l},$ $1 \leq k < l \leq
n$,
and defining relations:
\begin{align}
y_{k,l}^2 &= 1,\\
(y_{i,j} \,  y_{k, l})^2 &= 1,\\
(y_{i,k} \,  y_{k,j})^3 = (y_{i,j} \, y_{k,j})^3 &= (y_{i,k} \, y_{i,j})^3 = 1,
\end{align}
where distinct letters stand for distinct indices.
\end{thm}

\proof
We leave the proof to the reader, since one can  follow the same approach as in Theorems \ref{T:theorem1} and~\ref{T:flatH}.
The key point is that $GVB_n$ is the quotient of $FVB_n$ by the set of relations
\[
s_i \rr\ii = \rr\ii s_i, \quad i = 1, 2, \ldots, n-1.
\]
One can easily verify that it implies that~$y_{j,i} = y_{i,j}$, for~$1 \leq i < j \leq n$. Hence, $GH_n$ is generated by elements~$y_{k, l}$, for~$1 \leq k < l \leq n$. If we rewrite the set of relations of $FH_n$ in these generators and we proceed as in Proposition \ref{P:VirtualFlat} we get the set of relations
given in the statement. As before, one can also use the Reidemeister-Schreier method to check that this is a complete set of relations.

\endproof

As corollary, we have:

\begin{cor} 
The group $GVB_n$ is linear.
\end{cor}



\bibliography{Unrest-flatV3}{}
\bibliographystyle{plain}


\end{document}